\date{October 13, 2000} 
\title {     The Dimension of the Planar Brownian Frontier
is 4/3
}
\author {Gregory F.  Lawler\thanks{Duke University, Research
supported by the National Science Foundation}
\and Oded Schramm\thanks{Microsoft Research}
\and Wendelin Werner\thanks{Universit\'e Paris-Sud}
}
\documentclass[12pt]{article}
\usepackage{amsmath}
\usepackage{amsthm}
\usepackage{amsfonts}
\usepackage{graphicx}

\newif\ifdraft
\draftfalse
\numberwithin{equation}{section}

\newtheorem{theorem}{Theorem}
\numberwithin{theorem}{section}

\def\Prob{\Ps}
\newcommand{\R}{\mathbb{R}}
\newcommand{\C}{\mathbb{C}}

\def\H{\mathbb{H}}

\def\U{\Disk}

\def\Re{{\rm Re}}

\def \eps {\epsilon}
\def \tx {{\tilde \xi}}
\def \P {{\bf P}}
\def\Bb#1#2{{\def\md{\bigm| }#1\bigl[\,#2\,\bigr]}}
\def\BB#1#2{{\def\md{\Bigm| }#1\Bigl[\,#2\,\Bigr]}}
\def\Bs#1#2{{\def\md{\mid}#1[\,#2\,]}}
\def\Pb{\Bb\P}
\def\Eb{\Bb\E}

\def\EB{\BB\E}
\def\Ps{\Bs\P}

\def \E {{\bf E}}

\def\SG/{$SLE_\slepar$}
\def \SLE/{$SLE_6$}

\def \SLG/ {{\SG}}

\def \tx {{\tilde \xi}}

\long\def\hide#1{}

\def\ev#1{{\cal {#1}}} 

\def\slepar{\kappa}
\def\st{\,:\,}

\def\E{{\bf E}}

\def\tilde{\widetilde}
\def\mobtr/{M\"obius transformation}
\def\hat{\widehat}

\def \Disk {{\mathbb U}}

\begin{document}
\maketitle

\section{Introduction}

The purpose of this note is to announce and
sketch the proofs of
results determining the Hausdorff dimension of certain
subsets of planar Brownian paths.
Proofs are currently written down in a sequence
of preprints \cite{LSW1,LSW2,LSW3,LSWan}.
The present announcement will give an overview of some of the
contents of this series.
At this point, it appears that some parts of the proofs can
be significantly simplified using different means (see
\cite{Barcelona}).
However, it seems premature to give details about this
at this time.

Let $B_t$ be a
Brownian motion taking values in $\R^2$
(or $\C$). The {\em hull} of the Brownian motion $B$
at time $t$
is the union of $B[0,t]=\{B_s\st 0\le s\le t\}$
with the bounded components of the complement $\R^2\setminus B[0,t]$
of $B[0,t]$.
In other words, the hull is the Brownian path with all
the holes  ``filled in''.
The boundary of the hull is called the
{\em frontier} or {\em outer boundary} of
Brownian motion.  Based on
simulations and the analogy with self-avoiding
walks, Mandelbrot \cite{M} conjectured that the dimension
of the frontier is $4/3$.  In this note, we announce
a proof of this conjecture as well as proofs of some other
related questions on exceptional subsets of planar Brownian
paths.

A time $s \in (0,t)$ is called a {\em cut time}
and $B_s$ is called a {\em cut point} for $B[0,t]$
if
\[  B[0,s) \cap B(s,t] = \emptyset . \]
A point $B_s$ is called
a {\em pioneer point} if $B_s$ is on the frontier
at time $s$, that is, if $B_s$ is on the boundary
of the unbounded component of the complement of $B[0,s]$.

\begin{figure}\label{bm}
\centerline{\includegraphics*[height=2.5in]{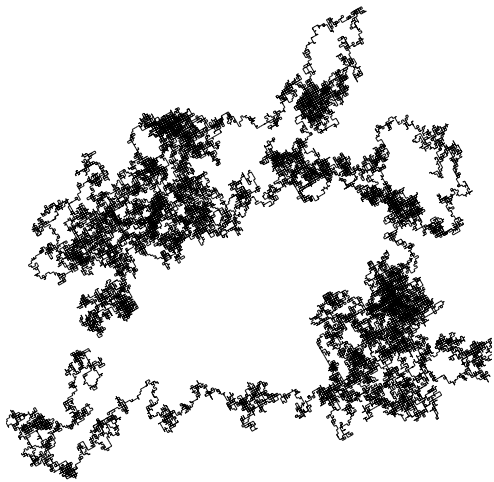}}
\vskip 0.5cm 
\caption{
{\bf Simulation of a planar Brownian path}}\bigskip
\end{figure}

\begin{theorem}[\cite{LSW2,LSWan}]\label{dimensions}
  Let $B_t$ be a planar Brownian motion.
With probability one,
the Hausdorff dimension of the frontier of $B[0,1]$
 is $4/3$; the
Hausdorff dimension of the set of cut points is $3/4$;
and the Hausdorff dimension of the set of pioneer
points is $7/4$.
\end{theorem}

This theorem is a corollary of a theorem determining
the values of the Brownian intersection exponents
$\xi(j,\lambda)$
which we define in the next section.  In \cite{Lfront,Lcut,Lbuda},
it had been established that the dimension of the
frontier, cut points, and pioneer points are $2 - \xi(2,0),
2 - \xi(1,1),$ and $2 - \xi(1,0)$, respectively.
Duplantier and Kwon \cite{DK} were the first to conjecture
the values $\xi(1,1) = 5/4, \xi(1,0) = 1/4$ using ideas
from conformal field theory.
Duplantier has also developed another
 non-rigorous approach to these
results based on ``quantum gravity'' (see e.g. \cite {Dqg}).
For a more complete list of references and background,
see e.g. \cite {LSW1}.

\section{Intersection exponents}

Let $B^1,B^2,\ldots,B^{j+k}$ be independent planar
Brownian motions with uniformly distributed starting points
on the unit circle.
Let $T^l_R$ denote the first time at which $B^l$ reaches
the circle of radius $R$, and let $\omega^j_R
  = B^j[0,T_R^j]$.
  The intersection exponent
$\xi(j,k)$ is defined by the relation 
\[   \Pb{(\omega_R^1 \cup \cdots \cup
   \omega_R^j) \cap (\omega_R^{j+1}
   \cup \cdots \cup\omega_R^{j+k})
   = \emptyset } \approx R^{-\xi(j,k)},\qquad R\to\infty
, \]
where $f\approx g$ means $\lim(\log f/\log g)=1$.
Using subadditivity, it is not hard to see that
there are constants $\xi(j,k)$ satisfying this relation.
Let
\[  Z_R = Z_R ( \omega^1_R, \ldots, \omega^j_R) :=
\Prob{(\omega_R^1 \cup \cdots \cup
   \omega_R^j) \cap \omega_R^{j+1} = \emptyset
   \md \omega_R^1,\ldots,\omega_R^j} . \]
Then 
\begin {equation}
\label {jl}
    \E[Z_R^k] \approx R^{-\xi(j,k)} ,\qquad R\to\infty. 
\end {equation}
In the latter formulation, there is no need to restrict
to integer $k$; this defines $\xi(j,\lambda)$ for
all $\lambda > 0$.  The existence of these exponents is
also very easy to establish. 
 The disconnection exponent
$\xi(j,0)$ is defined by the relation
\[   \Prob{Z_R > 0 } \approx R^{-\xi(j,0)},\qquad R\to\infty . \]
Note that $Z_R > 0$ if and only if 
$\omega_R^1 \cup \cdots \cup \omega_R^j$ 
do not disconnect the circle of radius 1 from 
the circle of radius $R$.
With this definition, $\xi(j,0) = \lim_{\lambda\searrow 0}\xi(j,\lambda)$
\cite{Lstrict}.

Another family of intersection exponents are the
half-plane exponents $\tilde{\xi}$. Let $\H$ denote
the open upper half-plane and define $\tilde{\xi}$ by
\[  \Prob{(\omega_R^1 \cup \cdots \cup
   \omega_R^j) \cap (\omega_R^{j+1}
   \cup \cdots \cup\omega_R^{j+k})
   = \emptyset ;\:  \omega_R^1\cup\cdots\cup\omega_R^{j+k}
   \subset \H } \approx R^{-\tx (j,k)} . \]

If $j_1,\ldots,j_n$ are positive integers and
$\lambda_0,\lambda_1,\ldots,\lambda_n \ge 0$, there is a
natural way to extend the above definitions and define the exponents
\[  \xi(j_1,\lambda_1,\ldots,j_n,\lambda_n),
\;\;\;\; \tilde \xi(\lambda_0,j_1,\ldots,j_n,\lambda_n). \]
In fact \cite{LW1}, there is a unique extension of the
intersection exponents \[ \tx(a_1,a_2,\ldots,a_n),\;\;\;\;
\xi(a_1,a_2,\ldots,a_n)\]
 to nonnegative reals $a_1,
\ldots,a_n$ (in the case of $\xi$, at least two of the
arguments must be at least $1$) such that: 
\begin {itemize}
\item
 the
exponents are symmetric functions
\item they satisfy
the ``cascade relations''
\begin{align}
\label {casc1}
   \tx(a_1,\ldots,a_{n+m})
&=
     \tx(a_1,a_2,\ldots, a_n,\tx(a_{n+1},\ldots,a_{n+m})) , \\
\label {casc2}
   \xi (a_1,\ldots,a_{n+m}) &=
     \xi(a_1,a_2,\ldots,a_n,\tx(a_{n+1},\ldots,a_{n+m})),\qquad
a_1,a_2\ge 1
 . \end{align}
\end {itemize}
Note that $\xi( j, \lambda)$ for $\lambda \in [0,1)$
and positive integer $j$ is defined directly via (\ref {jl}).

\begin{theorem}[\cite{LSW3}]
For all $a_1,a_2,\ldots, a_n \geq 0$,
\begin{equation}  \label{2}   \tx(a_1,\ldots,a_n) =
\frac {
\left(
\sqrt { 24 a_1 + 1} + \sqrt { 24 a_2 + 1 }
+ \cdots + \sqrt { 24 a_n +1 }
- (n-1) \right)^2  - 1  }{ 24 } .
\end{equation}
\end{theorem}

\begin{theorem}[\cite{LSW2,LSWan}]
For all $a_1,a_2,\ldots, a_n \geq 0$
with $a_1, a_2 \ge 1$,
\begin{equation}     \xi (a_1,\ldots,a_n) =
\frac {
\left(
\sqrt { 24 a_1 + 1} + \sqrt { 24 a_2 + 1 }
+ \cdots + \sqrt { 24 a_n +1 }
- n \right)^2  - 4  }{ 48 } .
\end{equation}
For all positive integers $j$ and all $\lambda \geq 0$,
\begin{equation}  \label{3}
  \xi(j,\lambda) = \frac{(\sqrt{24 j + 1} + \sqrt{24 \lambda
 + 1} - 2)^2 - 4}{48} .
\end{equation}
\end{theorem}

Note that three particular cases of the last theorem are
$\xi(2,0) = 2/3, \xi(1,1) = 5/4, \xi(1,0) = 1/4$, which, 
as noted above, are the
values from which Theorem~\ref{dimensions} follows.
The link between the exponents and the Hausdorff dimensions 
of the frontier, the  sets of cut-points 
and the set of pioneer points
loosely speaking 
goes as follows (here, for instance, for frontier points):
A point $x$ in the plane is in the $\epsilon$-neighborhood 
of a frontier  point if the Brownian motion (before time $1$)
reaches the $\epsilon$-neighborhood of $x$, and if the whole 
path $B[0,1]$ does not disconnect the disc of radius $\epsilon$
around $x$ from infinity.
This whole path can be divided in two parts: $B$ until it reaches
the circle of radius $\eps$ around $x$, and $B$ after it reaches this
circle. Both behave roughly like independent Brownian paths, and it is 
easy to see that 
the probability that $x$ is in the $\epsilon$ neighborhood of a 
frontier point decays like $\eps^{\xi(2,0)}$ when $\eps$ goes to 
zero.
This, together with some second moment estimates, implies that the
Hausdorff dimension of the frontier is $2- \xi (2,0)$.
See \cite {Lfront,Lcut,Lbuda} for details.

\section{Restriction property and universality}

This section reviews some of the results in \cite {LW2}.
Roughly speaking, a Brownian excursion in a domain is a Brownian motion
starting and ending on the boundary, conditioned to stay in
the domain.  To be more precise, consider the unit disk $D$,
start a Brownian motion uniformly on the circle of radius
$1-\epsilon$ and then stop the Brownian motion when
it reaches the unit circle.  This gives a probability
measure $\mu_\epsilon$ on paths.  The Brownian excursion
measure is the infinite measure
\[   \mu = \lim_{\epsilon \searrow 0} 2 \pi
   \epsilon^{-1} \mu_\epsilon .\]
If $A_1,A_2$ are closed disjoint arcs on the circle, then the
measure $\mu(A_1,A_2)$ of the set of paths starting at $A_1$ and ending at
$A_2$  is finite; moreover, as $A_1,A_2$ shrink,
$\mu(A_1,A_2)$ decays like $e^{-L(A_1,A_2)}$, where $L(A_1,A_2)=
L(A_1, A_2; D)$
denotes $\pi$ times the extremal distance between $A_1$ and $A_2$
(in the unit disk).
A particular Brownian excursion divides the unit disk into
two regions, say $U^+,U^-$, in such a way that the unit disk is
the union of $U^+,U^-$ and the ``hull'' of the excursion
(we choose $U^+$ and $U^-$ in such a way that 
the starting point of the excursion, $\partial U^+ \cap \partial D$,
the end-point of the excursion 
and $\partial U^- \cap \partial D$ are ordered clockwise on
the unit circle).

It can be shown \cite {LW2}
 that the Brownian excursion measure is
invariant under M\"obius transformations of the unit disk.
(Here we are considering two curves to be the same
if one can be obtained from the other by an increasing
reparameterization.)  Hence, by transporting via a conformal
map, the excursion measure
can be defined on any simply connected domain.
There is another important property
satisfied by this measure that we call the {\em restriction
property}.  Suppose $A_1,A_2$ are disjoint
arcs of the unit circle and $D'$ is a simply connected
 subset of $D$ whose
boundary includes $A_1,A_2$.  Consider the set
 of excursions $X(A_1,A_2,D')$
starting at $A_1$, ending at $A_2$, and staying in $D'$.
The  excursion measure 
gives two numbers associated with this set
of paths: $\mu(X(A_1,A_2,D'))$ and
$\mu(\phi(X(A_1,A_2,D')))$,
where $\phi$ is an arbitrary conformal
homeomorphism from $D'$ onto the unit disk.
The restriction property for the excursion measure states that these
two numbers are the same, for all such $D', A_1, A_2$.

Let $\ev X$ denote the space of all compact connected
subsets $X$ of the closed unit disk such that the intersection
of $X$ with the unit circle has two labeled connected components
$u^+$ and $u^-$ and the complement of $X$ in the plane is
connected.  The hull of a Brownian excursion
is an example for such an $X$.
Conformally invariant
measures on $\ev X$
that satisfy the restriction property and have a well-defined
crossing exponent (i.e., a number $\alpha$ such that the
measure of the set of paths from $A_1$ to $A_2$ decays
like $e^{-\alpha L(A_1,A_2)}$ as $A_1,A_2$ shrink to distinct points) are
called
{\em completely conformally invariant (CCI)}.  The
Brownian  excursion measure
$\mu$ gives a CCI measure with $\alpha= \alpha (\mu) = 1$. 
 Another CCI measure can be obtained by taking
$j$ independent excursions and 
considering the hull formed by their union;
this measure has crossing exponent $\alpha = j$.  Given a CCI
measure $\nu$, the intersection exponent $\tx(\lambda_1,
\nu,\lambda_2)$ can be defined by saying that as $A_1,A_2$ shrink
(i.e. as $L ( A_1, A_2) \to 0$)
\begin{align*}
 &\int \exp\Bigl(-\lambda_1 L(A_1,A_2;U^+) - \lambda_2 L(A_1,A_2;U^-)\Bigr)
\,d\nu(X) \\
&\qquad\qquad
\qquad\qquad
\approx
\exp\Bigl(-\tx(\lambda_1,\nu,\lambda_2) L(A_1,A_2)\Bigr).
\end{align*}
Here $U^+$ and $U^-$ denote the two components 
 of $X$ in the unit disk that have 
respectively 
$u^+$ and $u^-$ as parts 
of their boundary.
In the case of the Brownian excursion measure
$\mu$, $\alpha = 1$ and
 $\tx(\lambda_1,\mu,\lambda_2) = \tx(\lambda_1,
1,\lambda_2)$ where the latter denotes the Brownian intersection
exponent as in the previous section.  In \cite{LW2} it is
shown that any CCI measure $\nu$
with crossing exponent $\alpha$ acts
like the union of ``$\alpha$ Brownian excursions'' at least in
the sense that
for all $\lambda_1, \lambda_2 \ge 0$,
\[  \tx(\lambda_1,\nu,\lambda_2) = \tx(\lambda_1,\alpha,\lambda_2) .\]

The measure $\nu$ also induces 
naturally  a family of measures on excursions on annuli
bounded by the circles of radius $r$ and $1$.
For the Brownian excursion measure,
it corresponds 
to those Brownian excursions (in the unit disc)
stopped at their hitting time
of the circle of radius $r$
in case their reach it.
We can define a similar exponent $\xi(\nu,\lambda)$, and it can be
shown that $\xi(\nu,\lambda) = \xi(\alpha,\lambda)$
for all $\lambda \ge 1$.  In particular,
if the intersection exponents for one CCI measure could be determined,
then we would have the exponents for all CCI measures.

It is conjectured that the continuum limit of
percolation cluster boundaries  and self-avoiding walks
give CCI measures (or, at least, there are CCI measures in the
same universality class as these measures).  For example, the
well-known conjectures on self-avoiding walk exponents can be
reinterpreted as the conjecture that a self-avoiding walk
acts like $5/8$ Brownian motion.  (More precisely, two
nonintersecting self-avoiding walks act like two intersecting
Brownian paths.  This interpretation is similar to the
conjecture of Mandelbrot).  The restriction property
appears to be equivalent to the property of having ``central charge
zero'' from conformal field theory.

In terms of computing the Brownian intersection exponents, the
result in \cite{LW2} does not appear promising.  If one
could compute percolation or self-avoiding walk exponents, one
might be able to get the Brownian exponents; however, percolation
and self-avoiding walks exponents are even harder to compute
directly than Brownian intersection exponents. It is 
in fact not even rigorously known 
that they exist.  What was needed was
a conformally invariant process on paths for which one
can prove the restriction property and compute intersection
exponents. As we shall know point out,
 such a process exists --- the stochastic
L\"owner evolution at the parameter $\kappa = 6$.

\section{Stochastic L\"owner evolution}

The stochastic L\"owner evolution ($SLE_\kappa$) was first developed
in \cite{S} as a model for conformally invariant growth
($\kappa$ denotes a positive real).
There, two versions of $SLE_\kappa$ were introduced,
which we later called {\em radial}\/ and {\em chordal}\/ SLE.
Moreover, it was shown \cite{S} that radial
$SLE_2$ is the scaling limit of loop-erased
random walk (or Laplacian random walk),
if the scaling limit exists and is
conformally invariant, and it was conjectured that
chordal $SLE_6$ is the scaling limit of percolation cluster boundaries.
We now give the definition of chordal SLE, 
explain more precisely the conjecture relating chordal $SLE_6$
and percolation, and then define radial SLE.

The L\"owner equation relates a curve (or an increasing
family of sets) in
a complex domain to a continuous real-valued function.
SLE is obtained by choosing the real-valued function to
be a one-dimensional Brownian motion.

Let $\beta_t$ denote a standard one-dimensional Brownian
motion and let $W_t = \sqrt{\kappa} \beta_{ t}$ where
$\kappa > 0$.  Let $\H$ denote the upper half-plane,
and for $z \in \H$ consider the differential equation
\begin{equation}  \label{1}
    \partial_t g_t(z) = \frac{2}{g_t(z) - W_t}
\end{equation}
with initial condition $g_0(z) = z$.  The solution is
well-defined up to a (possibly infinite) time $T_z$ at
which $\lim_{t\nearrow T_z}g_{t}(z) - W_{T_z}=0$.  Let $D_t$ be the set of
$z$ such that $T_z > t$.  Then $g_t$ is the conformal
transformation of $D_t$ onto $\H$ with
\[       g_t(z) = z + \frac{2t}{z} + o(\frac{1}{z}),
   \;\;\;\;\;  z \rightarrow \infty . \]
The set $K_t = \H \setminus D_t$ is called the $SLE_\kappa$ hull.
  The
behavior of $SLE_{\kappa}$  depends strongly on
the value of $\kappa$.  One way to see this is to
fix a $z \in \H$ (or $z \in \partial \H$) and let
$Y_t =[ g_t(z) - W_t]/\sqrt{\kappa}$.  Then $Y_t$ satisfies the
stochastic differential equation
\[     dY_t = \frac{2}{\kappa Y_t } \; dt  - d \beta_t \]
defining the so-called 
 Bessel process of dimension $1 + 4 / \kappa$.
If  $\kappa < 4$, the Bessel process never hits the origin and 
 $  K_t$ is a simple curve; if $\kappa > 4$, every
point is eventually contained in the hull.  If $D$ is a
simply connected domain and $w,z$ are on the boundary, we
define $SLE_{\kappa}$ starting at $w$ and ending at $z$ by
 applying a conformal map from $\H$ onto $D$
which takes $0$ and $\infty$ to $w$ and $z$, respectively. 
(This assumes that $w$ and $z$ are ``nice'' boundary points;
otherwise, one has to consider prime ends instead.)
This defines a process uniquely
only up to increasing time change.  This is
{\em chordal SLE}.

\begin{figure}\label{percbound}
\centerline{\includegraphics*[height=2.5in]{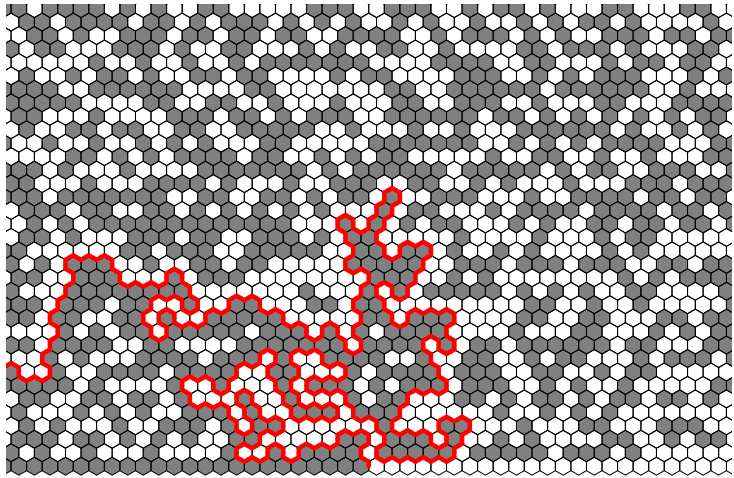}}
\vskip0.5cm
\caption{
{\bf The percolation boundary path}}\bigskip
\end{figure}

Consider critical percolation in $\H$ obtained
by taking each hexagon in a hexagonal grid of
mesh $\epsilon$ to be white or black with probability 
$1/2$, independently. This is in fact critical site
percolation on the triangular lattice.
Let $H_W$ be the connected component of the union
of the white hexagons and the positive real ray 
which contains the positive real ray, and let
$H_B$ be the connected component of the union of the
black hexagons and the negative real ray which
contains the negative real ray.
Then $\partial H_B\cap \partial H_W$ is a simple
path $\gamma_\epsilon$ connecting $0$ and $\infty$ in the closure of $\H$.
See figure~\ref{percbound}.

It is conjectured that when $\epsilon\searrow 0$ the law of the path
$\gamma_\epsilon$ tends to a conformally invariant measure.
Here, conformal invariance means that if we take an analogous
construction in another simply connected proper
subdomain of $\R^2$ in place of $\H$,
the resulting limit will be the same as the image of the limit
in $\H$ under a conformal homeomorphism between the domains.
Assuming this conjecture,
it can be shown \cite{Sinprep} that the scaling limit
of $\gamma_\epsilon$ is given by chordal $SLE_6$.
More precisely, the limit path, appropriately
parameterized, is given by
$\gamma(t)=g_t^{-1}(W_t)$,
where $g_t$ and $W_t$ are the $SLE_6$ maps and driving
parameter.

The restriction property for $SLE_6$ is an easy
consequence of the conjectured conformal invariance
and the independence properties of
critical percolation. However, one can actually prove \cite {LSW1}
 the
restriction property for $SLE_6$ without appealing
to the conjecture.
In contrast to the proof of the restriction
property for Brownian excursions, the proof of
the restriction property for $SLE_6$ is quite involved.
It is not hard to see that the restriction property
fails for $SLE_\kappa$ when $\kappa\ne 6$.

To understand the framework of the restriction
property proof, consider an infinitesimal
deformation of the domain (say, a small vertical slit $\alpha$)
away from the starting point of an $SLE$, and then
considers the conformal map taking the slit domain
to the entire half-space, the $SLE$ path is
transformed as well.  Essentially, it is mapped to a path
doing a L\"owner evolution with a different ``driving
function'' than $W_t$.   It is necessary to see how
this new driving function evolves as the slit $\alpha$ grows.
With the aid of stochastic calculus, it can be shown that
the new driving function is a martingale only
when $\kappa = 6$.  This argument is still rather mysterious.
An explicit calculation is done and the drift term disappears
exactly when $\kappa = 6$.

The definition of radial SLE is similar to chordal SLE.
Again, take $\beta_t$ to be a standard Brownian motion
on the real line.   Define $\zeta_t:= \exp(i\sqrt\kappa\beta_t)$,
which is a Brownian motion on the unit circle.
Take $g_0(z)=z$ for $z\in\U$,
and let $g_t(z)$ satisfy the L\"owner differential equation
$$
\partial_t g_t(z) =
g_t(z)\frac{g_t(z)+\zeta_t}{g_t(z)-\zeta_t}\,,
$$
up to the first time $T_z$ where $g_t(z)$ hits $\zeta_t$.
Let $D_t:=\{z\in\U\st T_z>t\}$ and $K_t:=\U\setminus D_t$.
This defines radial $SLE_\kappa$ from $1$ to $0$ in $\U$.
Note that the main difference from chordal SLE is
that $g_t$ is normalized at an interior point $0$ and the
set $K_t$ ``grows towards'' this interior point.

\section{Exponents for $SLE_6$}

Intersection exponents for $SLE_{\kappa}$ can
be computed for any $\kappa$.  We will consider only
the case $\kappa = 6$ here.
We will briefly describe how to compute the exponents 
$\tilde \xi (\lambda_1, SLE_6,\lambda_2)$ and 
$\xi (SLE_6, \lambda)$, but we first focus on another
family of exponents 
$\hat\xi$  that can
be described as follows.
Start a chordal \SLE/ at the top left corner $i\pi$
 of the rectangle $R:=[0,L] \times [0,\pi]$ 
going towards the lower right corner $L$.
Let $K$ be the hull of
the \SLE/ at the first time that
it hits the right edge $\{L\}\times[0,\pi]$.
Let ${\mathcal L}_+$ be $\pi$ times  the
extremal distance between the vertical edges of $R$
in $R\setminus K$, where we take ${\mathcal L}_+=\infty$
if there is no path in $R\setminus K$ connecting the
vertical edges of $R$ (that is, if $K$ intersects the lower edge).
Then the intersection exponent $\hat\xi(SLE_{6},\lambda)$
can be defined by the relation
\[         \E[e^{- \lambda {\cal L+}}] \approx e^{-\hat\xi(SLE_6
    ,\lambda) L},
    \;\;\;\;  L \rightarrow \infty . \]
The exponent $\hat\xi(SLE_6,\lambda)$ will 
turn out to be closely related
to the Brownian exponent $\tilde \xi (1/3, \lambda)$ and therefore 
also to
the Brownian exponent $\tx(1,\lambda)$.
However, observe that the SLE is here allowed to touch one horizontal
edge of the rectangle.  For this reason, we use the notation
$\hat \xi$, instead of $\tx$.

To calculate $\hat\xi(SLE_6,\lambda)$, we map the rectangle to the upper
half-plane
and consider the exponent there.  Let $\phi$ be the conformal
map taking $R$ onto $\H$ which satisfies
$\phi(L)=\infty$, $\phi(0)=1$, $\phi(L+i\pi)=0$,
and set $x:=\phi(i\pi)$.
Start a chordal \SLE/ in $\H$ from
$x$ to $\infty$.  Let $K = K_{T}$ be the hull of the
$SLE$ at the first time, $T=T_R$,
such that $K_t$ hits $(-\infty,0]\cup[1,\infty)$ (one can verify that
$T_R < \infty$).   Let ${\cal L} = {\cal L}_R$ be 
$\pi$ times the extremal distance between $(-\infty,0)$ and $(0,1]$
in $\H \setminus K$.  Then,
\[   \Eb{e^{-\lambda {\cal L}}} \approx e^{-L\hat\xi(SLE_6,\lambda)},
   \;\;\;\;  L \rightarrow \infty . \]
Let 
$$
f_t(z)=\frac{g_t(z)-g_t(0)}{g_t(1)-g_t(0)},
$$
i.e., $g_t$ normalized to fix $0$ and $1$.
It is not hard to verify that
$e^{\mathcal L}\approx 1-f_T(\sup(K_T\cap[0,1]))$.
In fact, we show in~\cite{LSW1} that
\begin {equation}
\label {calc}
\Eb{e^{-\lambda\mathcal L}}\approx\EB{\bigl((1-x) f_T'(1)\bigr)^\lambda}.
\end {equation}
Define
$$
w_t := \frac{W_t-g_t(0)}{g_t(1)-g_t(0)}, \ u_t= f_t' (1).
$$
The couple $(w_t, u_t)$ is solution of 
a stochastic differential equation, and it is easy to check that 
$T$ is in fact the first time at which $w_t$ hits 0 or 1.
Define $F(w,u) = {\bf E} [ u_T^\lambda | u_0 = u , w_0 = w]$.
The right hand side of (\ref {calc}) 
is then calculated via writing a PDE for $F$. 
The PDE becomes an ODE when symmetries are accounted for.
The solution of the ODE turns out to be a hypergeometric function
(the hypergeometric function predicted by 
Cardy \cite{Cardy} for the crossing probability of 
a rectangle for critical percolation is recovered as
the special case $\lambda=0$). 
This leads to
\[      \hat\xi(SLE_6,\lambda) = \frac{6 \lambda + 1 +
    \sqrt{24 \lambda +1}}{6} . \]
In particular, $\tx(SLE_6,0) = 1/3$ which can very
loosely
be interpreted as ``\SLE/ allowed to bounce
off of one side equals  a third of a Brownian
motion.''  
This leads directly, using the 
universality arguments developed in \cite {LW2}, to the 
fact that $\tx (1/3, \lambda) = \hat \xi (SLE_6, \lambda)$
for all $\lambda >0$.
In particular, $\tx (1/3, 1/3) = 1$, and the cascade relation then
gives $\tx (1, \lambda) = \hat \xi
 (SLE_6, \hat \xi (SLE_6, \lambda))$
and the value of these exponents
\cite {LSW1}.

By considering the regions above and below
$K$ in the rectangle $R$, and disallowing
$K$ to touch the horizontal sides, we can also 
define the two-sided exponents $\tx(\lambda_1,
SLE_6,\lambda_2)$.  As above, the calculation of this
exponent was translated to a question about a PDE.
However, in this case the PDE did not become an ODE,
and was not solved explicitly.  
But some eigenfunctions for the PDE were found,
which proved sufficient.  The leading eigenvalue
was shown to be equal to the sought exponent
\cite{LSW3}. The result is
\[  \tx(\lambda_1,SLE_6,\lambda_2) = 
\frac {
\left(
\sqrt { 24 \lambda_1 + 1} + 3
 + \sqrt { 24 \lambda_2 +1 }
 \right)^2  - 1  }{ 24 } . \]
The universality 
ideas show that  $\tx(\lambda_1,1,
\lambda_2) = \tx (\lambda_1, SLE_6, \lambda_2)$, so that  we 
get the value of $\tx (\lambda_1, 1, \lambda_2)$ 
for all $\lambda_1, \lambda_2 \ge 0$ and therefore, using the
cascade relations 
(\ref {casc1}),
(\ref {2}) follows.  

The exponent $\xi(SLE_6,\lambda)$ is defined by considering
the hull $K$ of
radial \SLE/ started on the unit circle stopped
at the first time it reaches the circle of radius
$r$ (when $r \to 0$).
  Let ${\cal L}$ be $\pi$ times the extremal distance
between the two circles in the annulus with $K$
removed.  Then, when $r \to 0$,
\[  \E[e^{-\lambda {\cal L}}] \approx r^{
   \xi(SLE_6,\lambda)} . \]
In a similar way the evaluation of this exponent
can be reduced to analyzing $\Eb{|g^\prime_t(z)|
^\lambda}$.  The calculation yields
\[  \xi(SLE_6,\lambda) =
   \frac{4\lambda +1 +\sqrt{24\lambda+1}}{8} ,
\]
at least for $\lambda\ge 1$.
From this and the ideas in \cite {LW2},
one can identify  $\xi (1, \lambda)$ 
with $\xi (SLE_6, \lambda)$ for all 
$\lambda \ge 1$. From this and the cascade relations
(\ref {casc2}),
(\ref{3}) for all $\lambda \ge 1$ 
and (\ref {2})  follow.

\section{Analyticity}

Unfortunately, the universality
argument from \cite{LW2} does not
show that $\xi(SLE_6,\lambda)=\xi(1,\lambda)$
when $\lambda<1$.
To derive the values of
$\xi(j,\lambda)$ for all $\lambda > 0$, we show
that the function $\lambda \mapsto \xi(j,\lambda)$ is
real analytic for $\lambda > 0$ \cite{LSWan}.
Then,
 the formula for $\xi(j,\lambda)$ where $\lambda<1$ follows by analytic
continuation.

To prove analyticity, we interpret $\xi(j,\lambda)$ as
an eigenvalue for an analytic function $\lambda
\mapsto T_\lambda$, taking values in the space
of bounded operators on a Banach
space.  For notational ease, let use assume
$j=1$.  Consider the set of continuous functions
$  \gamma:[0,t] \rightarrow \C $
with $\gamma(0) = 0 , |\gamma(t)| = 1$ and
$0 < |\gamma(s)| < 1$ for $0 < s < t$.
Let ${\cal C}$ be the subset of these functions with the
property that there is a connected component
$D$ of $\Disk \setminus \gamma$ whose boundary
includes both the origin and a subarc of the unit
circle (here we write $\gamma$ for $\gamma[0,t]$).
 Consider the following transformation on
$\gamma$.  If $\gamma \in {\cal C}$, start a Brownian
motion at the endpoint of $\gamma$ and let it run
until it hits the circle of radius $e$, producing
a curve $\beta$. Let  $\gamma_1$ be $\gamma\cup\beta$ rescaled
by $1/e$, so that its endpoint is once again on
the unit circle.  Let
$\Phi(\gamma_1)$ be the difference in the $\pi$-extremal
distance between $0$ and the unit circle
in $\Disk \setminus \gamma_1$ and the $\pi$-extremal
distance between $0$ and the circle $(1/e)\partial\U$ in
$\Disk \setminus (1/e) \gamma$.  (These distance are
both infinite but the difference can easily be defined.)
  If $f$ is a bounded function
on ${\cal C}$ and $z \in \C$ we define $T_zf$ by
\[   T_zf (\gamma) = \E[e^{-z \Phi(\gamma_1)}
           f(\gamma_1)] . \]
This is well-defined if $\Re(z) > 0$, and if $\lambda >0$,
then $e^{-\xi(1,\lambda)}$ is an eigenvalue for $T_\lambda$.
In fact, by appropriate choice of a Banach space norm, it
can be shown that for any $\lambda > 0$, there is a
complex neighborhood about $\lambda$ such that $z
\rightarrow T_z$ is holomorphic and such that $e^{-\xi(1,\lambda)}$
is an isolated point in the spectrum of $T_\lambda$.
The choice of the norm is similar to that used in
\cite{Ruelle} to show that the free energy of a
one-dimensional
Ising model with exponentially decaying interaction
is analytic in the temperature.  A coupling argument
is used in establishing exponential convergence that is
necessary to show that the eigenvalue is isolated.
Since the eigenvalue is isolated, it follows that it
depends analytically on $\lambda$.

\bigskip

\filbreak
\begingroup
\small
\parindent=0pt

\def\email#1{\texttt{#1}}
\vtop{
\hsize=2.3in
Greg Lawler\\
Department of Mathematics\\
Box 90320\\
Duke University\\
Durham NC 27708-0320, USA\\
\email{jose@math.duke.edu}
}
\bigskip
\vtop{
\hsize=2.3in
Oded Schramm\\
Microsoft Corporation,\\
One Microsoft Way,\\
Redmond, WA 98052; USA\\
\email{schramm@Microsoft.com}
}
\bigskip
\vtop{
\hsize=2.3in
Wendelin Werner\\
D\'epartement de Math\'ematiques\\
B\^at. 425\\
Universit\'e Paris-Sud\\
91405 ORSAY cedex, France\\
\email{wendelin.werner@math.u-psud.fr}
}
\endgroup

\filbreak

\end {document}